\documentclass[10pt,a4paper]{amsart}
\usepackage{amssymb,amsmath}
\usepackage{color}
\numberwithin{equation}{section}
\vspace{8cm}
\newtheorem{theorem}{Theorem}
\newtheorem{lemma}[theorem]{Lemma}
\newtheorem{prop}[theorem]{Proposition}

\newtheorem{oss}[theorem]{Remark}

\numberwithin{theorem}{section}
\newtheorem*{theorem*}{Theorem}
\newcommand{\Pk}{\mathcal{P}^+_k}
\newcommand{\Mpiu}{\mathcal{M}^+_{0,1}}

\newcommand{\R}{{\mathbb R}}

\newcommand{\refe}[1]{{\rm (\ref{#1})}}
\newcommand{\dis}{\displaystyle}

\newcommand{\e}{\'e\ }

\begin{document}
\parindent=0pt

\title[On the inequality  $F(x,D^2u)\geq f(u) +g(u)\, |Du|^q$] {On the inequality  $F(x,D^2u)\geq f(u) +g(u)\, |Du|^q$}
\author[I. Capuzzo Dolcetta, F. Leoni, A. Vitolo]
{Italo Capuzzo Dolcetta, Fabiana Leoni, Antonio Vitolo}

\address{ ICD and FL: Dipartimento di Matematica\newline
\indent Sapienza Universit\`a  di Roma}
\email{capuzzo@mat.uniroma1.it}
\email{leoni@mat.uniroma1.it}
\address{AV: Dipartimento di  Matematica\newline
Universit\`a di Salerno}
\email{vitolo@unisa.it}

\begin{abstract}
We consider  fully nonlinear degenerate elliptic equations with zero and first order terms. We provide a priori upper bounds and  characterize the existence of entire subsolutions under growth conditions on the lower order coefficients which extend the classical Keller--Osserman condition for semilinear equations.

 \end{abstract}\maketitle

\section{Introduction}\label{INTRO}

In 1957 J.B. Keller \cite{K} and R. Osserman \cite{O} simultaneously and independently proved that, for a given positive,  continuous and nondecreasing function $f$, the semilinear differential inequality
\begin{equation}\label{Isemi}
\Delta u \geq f(u)
\end{equation}
possesses an entire solution $u:\R^n\to \R$ if and only if 
\begin{equation}\label{IKO}
\int_0^{+\infty} \frac{dt}{\left(\int_0^t f(s)\, ds\right)^{1/2}} =+\infty\,\, \,.
\end{equation}
Moreover, if the Keller--Osserman condition \refe{IKO} fails and $u$ satisfies \refe{Isemi} in a proper open subset  $\Omega \subset \R^n$, then $u$ is bounded from above by a universal function of the distance from the boundary $\partial \Omega$, determined only by $f$.  
These results turned out to be rich of consequences and applications, and numerous generalizations have been established in the subsequent literature. Some of the most recent results will be recalled below.

Carrying on the study started in \cite{CDLV}, the aim of the present paper is to establish analogous results for the more general  fully nonlinear differential inequality
\begin{equation}\label{IeqF}
F(x, D^2u) \geq f(u) +g(u)\, |Du|^q\, ,\, x\in \R^n\,\, \,.
\end{equation}
Here, $f$ and $g$ are assumed to be continuous and monotone increasing, with $f$ positive, $q$ belongs to $(0,2]$ and $F$ is a second  order degenerate elliptic operator, that is a continuous function $F:\R^n \times \mathcal{S}_n \to \R$ satisfying $F(x,O)=0$ and the (normalized) ellipticity condition
$$
0\leq F(x,X+Y)-F(x,X)\leq \hbox{tr}(Y)\, ,\quad \forall\, x\in \R^n\, ,\ X\, ,\ Y\in \mathcal{S}_n\,,\ Y\geq O\, ,
$$
$\mathcal{S}_n$ being the space of symmetric $n\times n$ real matrices equipped with the usual partial ordering.

The model cases for $F$ that we have in mind are the degenerate maximal Pucci operator  $\Mpiu$ defined by
\begin{equation}\label{OP1}
\mathcal{M}^+_{0,1} (X) = \sum_{\mu_i>0} \mu_i(X)\, 
\end{equation}
and  
\begin{equation}\label{OP2}
\Pk (X) = \mu_{n-k+1}(X) +\ldots +\mu_n(X) 
\end{equation}

where $\mu_1(X) \leq \mu_2 (X) \leq \ldots \leq \mu_n(X)$ are the ordered eigenvalues of the matrix $X$.

The Pucci extremal operators,  in the uniformly elliptic case, have been extensively  studied  in a monograph by  L. Caffarelli and X. Cabr\e, see  \cite{CAFCA}. Let us recall here that the operator \refe{OP1} is maximal  in the class of  degenerate elliptic operators vanishing at $X=O$. In particular, for any $1\leq k\leq n$  and for all $X\in \mathcal{S}_n$, one has 
$$
\Pk (X)\leq \mathcal{M}^+_{0,1} (X)\,\,.
$$
Hence, if $u$ satisfies inequality 
\begin{equation}\label{IeqPk}
\Pk (D^2u)\geq  f(u) +g(u)\, |Du|^q\, ,
\end{equation}
or inequality \refe{IeqF}, then $u$ is also a solution of 
\begin{equation}\label{IeqM+}
\Mpiu (D^2u) \geq  f(u) +g(u)\, |Du|^q\,\,.
\end{equation}

As for the operators
$\Pk$,  we refer to the recent works of  R. Harvey and B. Lawson Jr \cite{HL} and L. Caffarelli, Y.Y. Li  and L. Nirenberg \cite{CLN},  see also M.E. Amendola, G. Galise, A. Vitolo \cite{AGV} and G. Galise \cite{G} and  the references therein for further results. We just point out here that  such degenerate operators arise in several frameworks, e.g. the  geometric problem of mean curvature evolution of manifolds with co-dimension greater than one, as in L. Ambrosio, H.M. Soner \cite{AS}, as well as in the PDE approach to the convex envelope problem,  see A. Oberman, L. Silvestre \cite{OS}.

Since the considered operators are in non-divergence form, a natural approach  to the analysis of  the partial differential inequality \refe{IeqF} is that of viscosity solutions. So,  by solution of \refe{IeqF} or \refe{IeqM+} or \refe{IeqPk} we always mean  an upper semicontinuous subsolution in the viscosity sense. 

Let us present now our results in a rather informal way. A key point is that while we assume $f$ to be positive we will not make a sign assumption on the function $g$, so that the first order terms can be of either "absorbing" or "reaction" type.

Let us discuss first the case where $\lim_{t\to +\infty} g(t) >0$. In this case, the necessary and sufficient "sublinearity" condition \refe{IKO} which rules the semilinear case \refe{Isemi} should be of course generalized in order to take in proper account the first order terms. We prove indeed that inequality \refe{IeqM+} possesses an entire viscosity solution (see Theorem \ref{entire} below) if and only if
\begin{equation}\label{IKO+}
q\leq 1\quad \hbox{and }\  \int_0^{+\infty}   \frac{dt}{\left( \int_0^t f(s)\, ds\right)^{1/2}+\left(\int_0^t g^+(s)\, ds\right)^{1/(2-q)}}=+\infty  \,\, .
\end{equation}
The same condition is also proved, in Theorem \ref{entirePk}, to be necessary and sufficient for the existence of an entire viscosity solution of \refe{IeqPk}, provided that $g(t)\geq 0$ for all $t\in \R$. In particular, we see that if $q>1$ then no entire subsolutions can exist, independently of how slow the growth of $f$ and $g$ is, whereas growth restrictions both on $f$ and $g$ are needed in the case  $q\leq 1$.
\\
Moreover, in Theorem \ref{upperestimate} we  show that if condition \refe{IKO+} is violated and  $u$ satisfies inequality \refe{IeqF} in a proper open subset $\Omega\subset \R^n$, then $u$ satisfies the universal pointwise upper bound
$$
u(x)\leq {\mathcal R}^{-1}(d(x))\, 
$$
where $d(x)={\rm dist}(x,\partial \Omega)$ and ${\mathcal R}$ is given (assuming for simplicity $g \geq 0$)  by the formula
$$
\mathcal{R} (a) =\left\{
\begin{array}{ll}
\dis 2\, \left(\frac{n}{2-q}\right)^{1/(2-q)} \int_a^{+\infty} \frac{dt}{  \left( \int_a^t f(s)\, ds\right)^{1/2}+\left( \int_a^t g (s)\, ds\right)^{1/(2-q)} } & \hbox{ if } q<2\\[6ex]
\dis  \sqrt{\frac{n}{2}} \int_a^{+\infty} \frac{dt}{ \left(  \int_a^t e^{\frac{2}{n}\int_s^t g (\tau)\, d\tau} f(s)\, ds\right)^{1/2}} & \hbox{ if } q=2
\end{array} \right.\, .
$$

Assume now that $\lim_{t\to +\infty} g(t)\leq 0$. In this case, the zero and the first order terms in inequality \refe{IeqM+} are competing with opposite signs. Our analysis shows that in this case entire viscosity subsolutions exist if and only if a relaxed version of  condition \refe{IKO} involving $f$, $g$ and $q$ holds true, namely
\begin{equation}\label{IKO0}
\int_0^{+\infty}  \frac{ dt}
 {\left( \int_0^t  e^{-2 \int_s^t \left( -\frac{g(\tau)}{f(\tau)}\right)^{2/q}f(\tau)\, d\tau} f(s)\, ds \right)^{1/2}}=+\infty\,\, .
\end{equation}
In particular, if $\lim_{t\to +\infty} g(t)< 0$, then \refe{IKO0} is proved to be equivalent to
\begin{equation}\label{IKO-}
\int_0^{+\infty}  \left[ \frac{ 1}
 {\left( \int_0^t  f(s)\, ds \right)^{1/2}}+ \frac{1}{f(t)^{1/q}}\right]\, dt =+\infty\,\,.\end{equation}
Note that the above condition becomes, for $q=2$, the "subquadratic" growth condition
$$
\int_0^{+\infty}   \frac{dt}{f(t)^{1/2}} =+\infty\,\, .
$$
Moreover, also in this case as in the previously discussed one, if condition \refe{IKO0} fails and $u$ satisfies \refe{IeqF} in any open subset $\Omega\subset \R^n$, then $u$ is universally estimated from above by an explicit function of the distance from the boundary determined by $f$, $g$ and $q$, see Theorem \ref{upperestimate}. 

As in the original papers \cite{K}, \cite{O} our strategy for proving the above results, see Theorem \ref{entire} and Theorem \ref{entirePk}, is based on comparison with radially symmetric solutions of \refe{IeqPk} and \refe{IeqM+}, after a detailed analysis of the existence of entire maximal solutions of the associated ordinary differential equation. As a matter of fact,   entire solutions of \refe{IeqM+} exist if and only if entire radially symmetric solutions exist. Remarkably, this fact is proved by a comparison argument which works also in the currently considered degenerate cases, without any a priori growth  assumptions at infinity on $u$.

Let us point out finally that, as a general fact, when the Keller--Osserman type conditions do not hold true and the zero order term $f(u)$ is an odd function satisfying the sign condition $f(u)\, u\geq 0$, one can obtain universal local bounds from above and from below. This property has been largely used in the literature to obtain existence results for entire solutions as well as existence of "large" solutions in bounded domains, that is solutions  blowing up on the boundary. \\ We did not investigate in this direction in the present paper and we just recall in this respect the result of 
H. Brezis \cite{B} about existence and uniqueness of entire solutions of  semilinear equations  with $f(u)=|u|^{p-1}u$, $p>1$. For subsequent extensions to more general divergence form principal parts and zero order terms we refer to L. Boccardo, T. Gallouet and J.L. Vazquez \cite{BGV1}, \cite{BGV2}, F. Leoni \cite{L}, F. Leoni, B. Pellacci \cite{LP} and  L. D'Ambrosio, E. Mitidieri \cite{DM}.\\ In the fully nonlinear framework, analogous results have been more recently obtained by M.J. Esteban, P. Felmer, A. Quaas \cite{EFQ}, G. Diaz \cite{D}  and G. Galise, A. Vitolo \cite{GV}, M. E. Amendola, G. Galise, A. Vitolo \cite{AGV2}, and for Hessian equations, involving the $k$-th elementary symmetric function of the eigenvalues $\mu_1(D^2u), \ldots ,\mu_n(D^2u)$ by  J. Bao, X. Ji \cite{BJ}, J. Bao, X. Ji, H. Li \cite{BJL}, Q. Jin, Y.Y. Li, H. Xu \cite{JLX}. For application to removable singularities results see also D. Labutin \cite{Denis}.

As far as equations with gradient terms are concerned, the analogous "absorbing" property of superlinear first order terms in semilinear elliptic equations was singled out first by J.M. Lasry and P.L. Lions \cite{LL} and then extensively studied, see e.g. S. Alarc\'on, J. Garc\'\i a--Meli\'an and A. Quaas \cite{AGQ} and  P. Felmer, A. Quaas and B. Sirakov \cite{FQS} for fully nonlinear uniformly elliptic equations with purely first oder terms of the form $h(|Du|)$.
Moreover,   some results obtained by A. Porretta \cite{P} on the existence of entire solutions of the equation
$$
-\Delta u +f(u)+g(u)\, |Du|^2= h(x)
$$
are closely related to the present paper. Let us observe that they  are obtained by performing the change of unknown $v=\int_u^{+\infty} e^{-\int_0^t g(s)\, ds}dt$ which can be used only when the power $q$ of the gradient term is 2.

\section{On the associated ODE}

In this section we perform a fairly complete qualitative analysis of the Cauchy problem\begin{equation}\label{Cpode}
\left\{\begin{array}{l}
\varphi''+\frac{c-1}r\,\varphi' =f(\varphi) +g(\varphi)\, |\varphi'|^q\\
\varphi(0)=a\\
\varphi'(0)=0
\end{array}\right.
\end{equation}where $f,g$ are continuous non decreasing functions and $c,q$ are positive real numbers.\\
As indicated in the Introduction this analysis is one of the basic tools in our approach to the study of entire solutions of the elliptic PDE
$$  F(x,D^2u)=f(u)+g(u)\, |Du|^q \,.$$
By solution of (\ref{Cpode}) in $[0,R)$ with $0 < R \le +\infty$, we mean here and in the sequel a function $\varphi \in C^2((0,R)) \cap C([0,R))$ satisfying, moreover,
 $$0=\varphi' (0)=\lim_{r\to 0} \varphi'(r)\quad,\quad \exists\, \lim_{r \to 0^+}\varphi''(r)\in \R\,.$$ 
 Therefore the  ordinary differential equation in (\ref{Cpode}) has to be satisfied for $r=0$, too.  \\
Let us observe that the existence of local solutions of (\ref{Cpode})  follows from the standard theory of ordinary differential equations with continuous data.
\medskip

\begin{lemma}\label{convex-increasing} 
Let  $c>0$, $q>0$, and $f$ and $g$  be continuous. Let  $\varphi $ be a solution of \refe{Cpode} in $[0,R)$; then:
\begin{itemize}
\item[(i)] if $f$  is positive, then $\varphi$ is strictly increasing;
\item[(ii)] if $f$  is positive and $f,\ g$ are non decreasing,  and $c\geq 1$, then $\varphi$ is convex and
\begin{equation}\label{upg-}
\varphi'(r) \leq \left( \frac{f(\varphi (r))}{g^- (\varphi (r))}\right)^{1/q}\, \mbox{for all}\; r\in [0,R), 
\end{equation}
\item[(iii)] if $f$  is positive,   $f,\ g$ are  non decreasing, $g$ is non negative  and $c\geq 1$, then  
\begin{equation}\label{sc}
\varphi''(r)\geq \frac{\varphi '(r)}{r}\, \, \mbox{for all}\; r\in [0,R)\;.
\end{equation}
\end{itemize}
\end{lemma}
\noindent{\bf Proof.}  From \refe{Cpode} it follows that
$$
c\varphi''(0)= \lim_{r\to 0^+}\left( \varphi''(r)+(c-1)\,\frac{\varphi'(r)}r\right)=f(a)>0\, ,
$$ 
so that $\varphi'$ is increasing, hence positive,  in some interval $(0,r_0)$. Actually, one has $\varphi'(r)>0$ in the whole interval $(0,R)$, since, if not, there should be a point $r^*\in (0,R)$ satisfying  $\varphi'(r^*)=0$ and $\varphi''(r^*)\leq 0$ on the one hand, and $\varphi''(r^*)=f(\varphi(r^*))>0$ on the other hand. Hence  (i)  is proved.\\
Next, let us  prove (ii). Since $\varphi''(0)>0$,  there exists some $r_1>0$ such that $\varphi'(r)>0$ and $\varphi''(r)>0$  for $r\in (0,r_1]$. By contradiction, let us assume that there exists $\tau >r_1$ such that $\varphi''(\tau)<0$. Then,  the function $\varphi'$ has a strict local maximum point $r_0\in (0,\tau)$ and the set $\mathcal{R}=\{ r\in (0,r_0)\, :\, \varphi'(r)=\varphi'(\tau)\}$ is non empty. Let $\sigma =\min \mathcal{R}$, so that $\varphi'' (\sigma)\geq 0$. Therefore, we have found $\sigma\, ,\ \tau \in (0,R)$ such that
\begin{equation}\label{fipfis}
\sigma<\tau\, ,\qquad \varphi'(\sigma)=\varphi'(\tau)\, ,\qquad \varphi''(\sigma)-\varphi''(\tau)>0\, .
\end{equation}
On the other hand, equation \refe{Cpode} tested at $\sigma$ and $\tau$ yields
$$
\begin{array}{ll}
\varphi'' (\sigma)-\varphi''(\tau)=  & (c-1)  \left( \frac{\varphi'(\tau)}{\tau}-\frac{\varphi'(\sigma)}{\sigma}\right) \\[2ex]
& +\,  g(\varphi(\sigma))\varphi' (\sigma)^q -g(\varphi(\tau))\varphi' (\tau)^q \\[2ex]
& + \, f(\varphi (\sigma))-f(\varphi (\tau))\, ,
\end{array}
$$
In view of the monotonicity of $f$, $g$ and $\varphi$ and the assumption $c\geq 1$, this contradicts \refe{fipfis}. Therefore, $\varphi$ is convex in $[0,R)$ and from equation \refe{Cpode} we immediately obtain
$$
f(\varphi )+ g(\varphi ) (\varphi ')^q\geq 0 \quad \hbox{ in } [0,R)\, ,
$$
which yields \refe{upg-}. This proves (ii).\\
Finally, let us  prove (iii).  Multiplying equation \refe{Cpode} by $r^{c-1}$ and integrating between $0$ and $r$ yields

$$
r^{c-1}\varphi'(r)\leq \left[ f(\varphi(r))+g(\varphi (r))\varphi'(r)^q\right]\int_0^r s^{c-1}ds=\left[ \varphi'' (r) +\frac{(c-1)}{r}\varphi'(r)\right] \frac{r^c}{c}\, ,
$$
since $\varphi$, $\varphi'$, $f$ and $g$ are non decreasing.
Hence \refe{sc} is proved.

\hfill$\Box$ 

For the next results we focus on the case in which the data of problem \refe{Cpode} satisfy 
\begin{equation} \label{ipotesi}
q\in (0,2], c\geq 1, f, g \; \mbox{continuous and non decreasing},\; f \;\mbox{positive}
\end{equation}
and we obtain sharp estimates from above and from below on the function $\varphi'$.

\begin{lemma}\label{stimefip}
Assume (\ref {ipotesi}). If $\varphi$ is a solution of \refe{Cpode} in $[0,R)$, then, for every $r\in [0,R)$, 
\begin{equation}\label{ubg+}
\varphi'(r) \leq \left\{
\begin{array}{ll}
\dis 2^{2/(2-q)} \left[ \left(\int_a^{\varphi(r)} f(t)\, dt\right)^{1/2} + \left( \int_a^{\varphi(r)}g^+(t)\, dt\right)^{1/(2-q)}\right] & \hbox{ if } q<2\\[4ex]
\dis  \left( 2 \int_a^{\varphi(r)} e^{2\int_t^{\varphi(r)} g^+(s)\, ds} f(t)\, dt\right)^{1/2} & \hbox{ if } q=2
\end{array} \right.
\end{equation}
and also, if $g(a)\geq 0$,
\begin{equation}\label{lbg+}
\varphi'(r) \geq \left\{
\begin{array}{ll}
\dis \frac{1}{2} \left(\frac{2-q}{c}\right)^{1/(2-q)} \left[ \left( \int_a^{\varphi(r)}f(t)\, dt\right)^{1/2}+ \left( \int_{a}^{\varphi(r)} g^+(t)\, dt\right)^{1/(2-q)}\right] & \hbox{ if } q<2\\[4ex]
\dis \left( \frac{2}{c} \int_{a}^{\varphi(r)} e^{\frac{2}{c} \int_t^{\varphi(r)} g^+(s)\, ds} f(t)\, dt\right)^{1/2}& \hbox{ if } q=2
\end{array} \right.
\end{equation}

Furthermore, 
\begin{equation}\label{lbg-}
\varphi'(r) \geq \left( \frac{q}{c^{2/q}} \int_a^{\varphi(r)}  e^{-2 \int_t^{\varphi(r)} \left( \frac{g^-(s)}{f(s)}\right)^{2/q}f(s)\, ds}f(t)\, dt \right)^{1/2}
\end{equation}
and also, if $g(t)\leq 0$ for all $t\in \R$,
\begin{equation}\label{ubg-}
\varphi'(r) \leq  \left(2 \int_a^{\varphi(r)}   e^{-2 \int_t^{\varphi(r)} \left( \frac{g^-(s)}{f(s)}\right)^{2/q}f(s)\, ds}f(t)\, dt \right)^{1/2}
\end{equation}
Finally, if $\varphi $ is a maximal solution of \refe{Cpode} on $[0,R)$, with $R\leq +\infty$, then 
\begin{equation}\label{blowupfi}
\lim_{r \to R^-}\varphi(r)= +\infty\, .
\end{equation}
\end{lemma}

\noindent {\bf Proof.} By Lemma \ref{convex-increasing}, $\varphi$ is convex and increasing in $[0,R)$, so that, from the equation in \refe{Cpode} we immediately deduce
\begin{equation}\label{1}
\varphi ''  +g^- (\varphi) (\varphi')^q \leq  f(\varphi)+g^+(\varphi) (\varphi')^q\, .
\end{equation}
In particular, one has
\begin{equation}\label{2}
\varphi ''   \leq  f(\varphi)+g^+(\varphi) (\varphi')^q\, .
\end{equation}
Consider first the case $q<2$. Multiplying \refe{2} by $\varphi'$ and integrating in $(0,r)$, jointly with the increasing monotonicity of $\varphi'$, then yields
$$
\frac{(\varphi'(r))^2}{2}\leq  \int_a^{\varphi(r)} f(t)\, dt +(\varphi'(r))^q \int_a^{\varphi (r)} g^+(t)\, dt\, .
$$
By Young inequality with exponent $2/q>1$, we then obtain
$$
\frac{(\varphi'(r))^2}{2}\leq  \int_a^{\varphi(r)} f(t)\, dt +\frac{(\varphi'(r))^2}{4}+ 4^{q/(2-q)} \left( \int_a^{\varphi (r)} g^+(t)\, dt\right)^{2/(2-q)}\, ,
$$
which immediately gives \refe{ubg+} in the case $q<2$.
\newline If $q=2$, we multiply \refe{2} by $2 e^{-2 \int_a^{\varphi (r)} g^+(t)\, dt} \varphi'(r)$ and deduce
$$
\left( e^{-2 \int_a^{\varphi (r)} g^+(t)\, dt} (\varphi'(r))^2\right)' \leq 2 \,  e^{-2 \int_a^{\varphi (r)} g^+(t)\, dt} \varphi' (r) f(\varphi(r)) \, .
$$
Integrating in $(0,r)$ then yields \refe{ubg+} in the case $q=2$.

On the other hand,   the equation in \refe{Cpode} may be written as
$$
\frac{\left( r^{c-1}\varphi'(r)\right)'}{r^{c-1}} =f(\varphi(r))+ g(\varphi(r)) \varphi'(r)^q\, ,
$$
hence
$$
\left( r^{c-1}\varphi'(r)\right)'\leq \left[ f(\varphi(r))+ g^+(\varphi(r)) \varphi'(r)^q\right]\, r^{c-1}\, ,
$$
and  integration  in $(0,r)$, by the monotonicity of $\varphi$, $\varphi'$, $f$ and $g^+$, yields
$$
\frac{\varphi'(r)}{r} \leq \frac{1}{c} \left[  f(\varphi(r)) + g^+(\varphi(r)) (\varphi'(r))^q\right] \, .
$$
The above inequality  inserted in \refe{Cpode} then implies
\begin{equation}\label{lb}
\varphi''(r) +g^-(\varphi(r))(\varphi'(r))^q \geq \frac{1}{c} \left[ f(\varphi(r))+ g^+(\varphi(r)) (\varphi'(r))^q\right] \, .
\end{equation}
Assume  now that  $g(a)\geq 0$; in this case $g(\varphi (r))\geq 0$ for every $r\in [0,R)$ and the above inequality becomes
\begin{equation}\label{lb1}
\varphi''(r) \geq \frac{1}{c} \left[ f(\varphi(r))+ g^+(\varphi(r)) (\varphi'(r))^q\right] \, .
\end{equation}
If $q=2$, inequality \refe{lb1} can be easily integrated as before yielding  inequality \refe{lbg+} for $q=2$.

If $q<2$, \refe{lb1} can be split into the two 
$$
\varphi''(r)\geq \frac{f(\varphi (r))}{c}\, ,\qquad \varphi''(r) \geq \frac{g^+(\varphi(r)) (\varphi'(r))^q}{c}  \, .
$$
From the former  it easily follows that
$$
\varphi'(r) \geq \left(\frac{2}{c} \int_a^{\varphi(r)} f(t)\, dt\right)^{1/2}\, ,
$$
whereas the integration of the latter yields
$$
\varphi'(r)\geq 
 \left( \frac{2-q}{c} \int_a^{\varphi(r)} g^+(t)\, dt\right)^{1/(2-q)}\, .
$$
Adding term by term we obtain  \refe{lbg+} also for $q<2$. \\

Furthermore,  we notice that \refe{lb}  also gives
\begin{equation}\label{lb2}
\varphi''(r) +g^-(\varphi(r))(\varphi'(r))^q \geq \frac{f(\varphi(r))}{c}\, .
\end{equation}
If $q<2$, by Young inequality with exponent $2/q>1$ we obtain
$$
\varphi''(r) +\frac{q}{2} c^{2/q -1} \frac{g^-(\varphi)^{2/q}}{f(\varphi)^{2/q-1}}(\varphi') ^2+\left( 1-\frac{q}{2}\right) \frac{f(\varphi)}{c} \geq \frac{f(\varphi)}{c}\, ;
$$
hence
$$
\varphi''(r) +  \frac{g^-(\varphi)^{2/q}}{f(\varphi)^{2/q-1}}(\varphi') ^2 \geq \frac{q}{2} \frac{f(\varphi)}{c^{2/q}}\, ,
$$
which holds true, by \refe{lb2}, also for $q=2$. \\
By multiplying both sides by $2 \varphi' e^{2 \int_a^{\varphi} \left( \frac{g^-(t)}{f(t)}\right)^{2/q}f(t)\, dt}$ and by integrating in $(0,r)$ we obtain
$$
e^{2 \int_a^{\varphi(r)} \left( \frac{g^-(t)}{f(t)}\right)^{2/q}f(t)\, dt} (\varphi')^2 \geq \frac{q}{c^{2/q}}  \int_a^{\varphi(r)} e^{2 \int_a^{t} \left( \frac{g^-(s)}{f(s)}\right)^{2/q}f(s)\, ds} f(t)\, dt\, ,
$$
that is \refe{lbg-}.

On the other hand, if $g$ is nonpositive, inequality \refe{1} reads
$$
\varphi'' + \frac{g^-(\varphi)}{(\varphi')^{2-q}} (\varphi')^2 \leq f(\varphi)\, ,
$$
which, on the account of \refe{upg-}, implies
$$
\varphi'' + \frac{g^-(\varphi)^{2/q}}{f(\varphi)^{2/q-1}} (\varphi')^2 \leq f(\varphi)\, .
$$
Multiplication  by $2 \varphi' e^{2 \int_a^{\varphi} \left( \frac{g^-(t)}{f(t)}\right)^{2/q}f(t)\, dt}$ and  integration in $(0,r)$ yields, as before,
$$
e^{2 \int_a^{\varphi(r)} \left( \frac{g^- (t)}{f(t)}\right)^{2/q}f(t)\, dt} (\varphi')^2 \leq 2  \int_a^{\varphi(r)} e^{2 \int_a^{t} \left( \frac{g^- (s)}{f(s)}\right)^{2/q}f(s)\, ds} f(t)\, dt\, ,
$$
and \refe{ubg-} is proved.

Finally, let us prove \refe{blowupfi}. Since $\varphi$ is  convex and  increasing in $[0,R)$,   the limits $\lim_{r\to R^-} \varphi(r)$ and $\lim_{r\to R^-} \varphi'(r)$ exist. Now, if $R=+\infty$, then \refe{blowupfi} follows by convexity and strictly increasing monotonicity. On the other hand,  from estimate \refe{ubg+} it follows that if $\varphi$ is bounded in $[0,R)$, then  $\varphi'$ is bounded as well, and this contradicts the maximality of $R$, if $R<+\infty$\,. Hence, \refe{blowupfi} holds true in any case.
\hfill$\Box$

For the sequel, it is convenient to rewrite the estimates of Lemma \ref{stimefip} in the following equivalent formulations.

\begin{lemma}\label{stimer}
Assume \refe{ipotesi}. If $\varphi$ is a solution of \refe{Cpode} for some $a\in \R$ and $R>0$, then for every $r\in [0,R)$ one has
\begin{equation}\label{lbrg+}
r \geq \left\{
\begin{array}{ll}
\dis   \frac{1}{2^{2/(2-q)}}  \int_a^{\varphi (r)} \frac{dt}{ \left[ \left(\int_a^t f(s)\, ds\right)^{1/2} + \left( \int_a^t g^+(s)\, ds\right)^{1/(2-q)}\right] }  & \hbox{ if } q<2\\[4ex]
\dis  \int_a^{\varphi(r)} \frac{dt} {\left( 2 \int_a^t e^{2\int_s^t g^+(r)\, dr}f(s)\, ds\right)^{1/2}}   & \hbox{ if } q=2
\end{array} \right.
\end{equation}
as well as, if $g(a)\geq 0$, 
\begin{equation}\label{ubrg+}
r \leq  \left\{
\begin{array}{ll}
\dis 2\, \left(\frac{c}{2-q}\right)^{1/(2-q)} \int_a^{\varphi(r)} \frac{dt}{  \left( \int_a^t f(s)\, ds\right)^{1/2}+\left( \int_a^t g^+ (s)\, ds\right)^{1/(2-q)} } & \hbox{ if } q<2\\[6ex]
\dis  \sqrt{\frac{c}{2}} \int_a^{\varphi(r)} \frac{dt}{ \left(  \int_a^t e^{\frac{2}{c}\int_s^t g^+ (\tau)\, d\tau} f(s)\, ds\right)^{1/2}} & \hbox{ if } q=2
\end{array} \right.\, .
\end{equation}
Furthermore, for any $q\in(0,2]$ and $r\in (0,R)$, one has
\begin{equation}\label{ubrg-}
 r\leq \frac{c^{1/q}}{\sqrt{q}} \int_a^{\varphi(r)}   \frac{ dt}
 {\left( \int_a^t  e^{-2 \int_s^t \left( \frac{g^-(\tau)}{f(\tau)}\right)^{2/q}f(\tau)\, d\tau} f(s)\, ds \right)^{1/2}} \, ,
\end{equation}
and, if $g(t)\leq 0$ for all $t\in \R$,
\begin{equation}\label{lbrg-}
 r\geq \frac{1}{\sqrt{2}} \int_a^{\varphi(r)}   \frac{ dt}
 {\left( \int_a^t  e^{-2 \int_s^t \left( \frac{g^- (\tau)}{f(\tau)}\right)^{2/q}f(\tau)\, d\tau} f(s)\, ds \right)^{1/2}}\, .
\end{equation}
\end{lemma}
\noindent {\bf Proof.} We rewrite  estimate \refe{ubg+} in the form
$$
1\geq \left\{
\begin{array}{ll}
\dis \frac{\varphi'(r)}{ 2^{2/(2-q)} \left[\left( \int_a^{\varphi(r)} f(t)\, dt \right)^{1/2}+ \left( \int_a^{\varphi(r)}g^+(t)\, dt\right)^{1/(2-q)}\right] } & \hbox{ if } q<2\\[4ex]
\dis \frac{\varphi'(r)}{ \left(2  \int_a^{\varphi(r)} e^{2 \int_t^{\varphi(r)} g^+(s)\, ds}f(t)\, dt\right)^{1/2}} & \hbox{ if } q=2
\end{array} \right.
$$
Hence, integration  in $[0,r]$ yields \refe{lbrg+}.
The remaining statements \refe{ubrg+}, \refe{ubrg-} and \refe{lbrg-}  are deduced in a completely analogous way  from \refe{lbg+}, \refe{lbg-}  and  \refe{ubg-} respectively.

\hfill$\Box$

The next result Theorem \ref{R} identifies in the general framework of assumption \refe{ipotesi},  some integrability conditions on the data which are necessary and sufficient for the existence of maximal solutions of  the Cauchy problem \refe{Cpode} defined on the whole interval $[0,+\infty)$.\\
 Its proof makes use in particular of the following elementary calculus result:
\begin{lemma}\label{R(a)}
Let  $f$ and $h$ be continuous functions in $\R$, with $f$ positive and non decreasing, and $h$ nonnegative and non increasing. For every $a\in \R$ let us set
$$
{\mathcal R}(a) =\int_a^{+\infty} \frac{dt}{\left( \int_a^t e^{-2 \int_s^t h(r)\, dr}f(s)\, ds\right)^{1/2}}\, .
$$
Then, either ${\mathcal R(a)} \equiv +\infty$, or ${\mathcal R}(a)<+\infty$ for every $a\in \R$. In the latter case, ${\mathcal R} :\R \to (0,+\infty)$ is monotone non increasing and 
\begin{equation}\label{infy}
\lim_{a\to +\infty} {\mathcal R}(a)=0\, .
\end{equation}
\end{lemma}
The proof of this Lemma will be detailed below for the convenience of the reader.

\begin{theorem}\label{R}
Under assumption \refe{ipotesi},  let $\varphi$ be a maximal solution of the initial value problem \refe{Cpode}. Then,
\begin{itemize}

\item[(i)] if  $\dis \lim_{t\to +\infty} g(t)> 0$, then $\varphi$ is globally defined  in $[0,+\infty )$ if and only if 
\begin{equation}\label{ns+}
 q\leq 1 \quad \hbox{ and }\ \int_0^{+\infty}   \frac{dt}{(t\, f(t) )^{1/2}+(t\, g^+(t) )^{1/(2-q)}}=+\infty   
\end{equation} 

\item[(ii)] 
 if $\dis \lim_{t\to +\infty} g(t)\leq  0$, then $\varphi$ is globally defined  in $[0,+\infty )$ if and only if 
 \begin{equation}\label{ns0}
\int_0^{+\infty}  \frac{ dt}
 {\left( \int_0^t  e^{-2 \int_s^t \left( \frac{g^-(\tau)}{f(\tau)}\right)^{2/q}f(\tau)\, d\tau} f(s)\, ds \right)^{1/2}}\,  dt=+\infty\, 
\end{equation}
\item[(iii)]  if  $\dis \lim_{t\to +\infty} g(t)< 0$, then \refe{ns0} is equivalent to
\begin{equation}\label{ns-}
\int_0^{+\infty} \left[  \frac{1}{(t\, f(t))^{1/2}} +   \frac{1}{f(t)^{1/q}}\right]\, dt =+\infty\, 
\end{equation}

\end{itemize}
\end{theorem}
\noindent{\bf Proof.} (i) Assume $\lim_{t\to +\infty} g(t)>0$. Let us first observe that condition \refe{ns+}  is equivalent to
$$
 q\leq 1 \quad \hbox{ and }\  \int_{t_0}^{+\infty} \frac{dt}{(t\, f(t))^{1/2}+(t\, g^+(t))^{1/(2-q)}} =+\infty  \  \hbox{ for all } t_0\geq 0\, .
$$
By observing that, for all $a\in \R$ and $t\geq 2|a|$ one has
$$
\frac{t}{2}\, f\left( \frac{t}{2}\right) \leq \int_a^t f(s)\, ds\leq f(t)\, (t-a)\leq 2\, t\, f(t)\, ,
$$
as well as
$$
\frac{t}{2}\, g^+\left( \frac{t}{2}\right) \leq \int_a^t g^+(s)\, ds\leq g^+(t)\, (t-a)\leq 2\, t\, g^+(t)\, ,
$$
the above condition is then equivalent to
$$
 q\leq 1 \quad \hbox{ and }\  \int_{2|a|}^{+\infty} \frac{dt}{\left( \int_a^t f(s)\, ds\right)^{1/2}+ \left( \int_a^t g^+(s)\, ds\right)^{1/(2-q)}} =+\infty\  \hbox{ for all } a\in \R\,,
$$
and then also to 
\begin{equation}\label{ns+a}
 q\leq 1 \quad \hbox{ and }\  \int_a^{+\infty} \frac{dt}{\left( \int_a^t f(s)\, ds\right)^{1/2}+ \left( \int_a^t g^+(s)\, ds\right)^{1/(2-q)}} =+\infty\  \hbox{ for all } a\in \R\, .
\end{equation}
 Now, let  \refe{ns+}, and therefore \refe{ns+a}, be satisfied. If $\varphi\in C^2([0,R))$ is a maximal solution of \refe{Cpode}  for some $a\in \R$, then, by letting $r\to R^-$ in  \refe{lbrg+} with $q\leq 1$, and by using  \refe{blowupfi} and \refe{ns+a}  we immediately infer $R=+\infty$.

Conversely, assume that for any $a\in \R$ any maximal solution of the Cauchy problem \refe{Cpode} belongs to $C^2([0,+\infty))$, and let us select a maximal solution $\varphi$ satisfying the initial condition $\varphi(0)=a>0$ such that $g(a)> 0$. By inequality \refe{ubrg+}  it follows that, if $q>1$, then
$$
r \leq  \left\{
\begin{array}{ll}
\dis 2 \left( \frac{c}{2-q}\right)^{1/(2-q)} \int_a^{+\infty} \frac{dt}{\left( f(a) (t-a)\right)^{1/2} +\left( g^+(a) (t-a)\right)^{1/(2-q)}} & \hbox{ if } 1< q<2
\\[4ex]
\dis\sqrt{\frac{g^+(a)}{f(a)}}  \int_a^{+\infty} \frac{dt}{\sqrt{e^{\frac{2}{c} g^+(a) (t-a)}-1}} & \hbox{ if } q=2
\end{array} \right.
$$
which is a contradiction to the unboundedness of $r$. Hence, one has $q\leq 1$, and letting $r\to +\infty$ in \refe{ubrg+}, we obtain, by \refe{blowupfi}, that \refe{ns+a} is satisfied for sufficiently large $a$, and, therefore, \refe{ns+} holds true.

(ii) Assume that $\lim_{t\to +\infty} g(t)\leq 0$, so that $g(t)\leq 0$ for all $t\in \R$, and  let $\varphi \in C^2([0,R))$ be any maximal solution of the Cauchy problem \refe{Cpode}. By inequalities \refe{ubrg-} and \refe{lbrg-},  and by \refe{blowupfi},  we immediately infer that $R=+\infty$ if and only if, for all $a\in \R$,
$$
\int_a^{+\infty} \frac{dt}{\left( \int_a^t e^{-2\int_s^t \left(\frac{g^-(r)}{f(r)}\right)^{2/q}f(r)\, dr} f(s)\, ds\right)^{1/2}} =+\infty\, .
$$
By applying Lemma \ref{R(a)}  with $h(r)=\left( g^-(r)/f(r)\right)^{2/q}f(r)$, we deduce that  the above condition is equivalent to \refe{ns0}.

(iii) It remains to prove that, if $\lim_{t\to +\infty} g(t)=g_\infty<0$,  then \refe{ns0} is equivalent to \refe{ns-}. First of all, let us rewrite condition \refe{ns0} as
\begin{equation}\label{nsh}
\int_a^{+\infty} \frac{dt}{\left(\int_a^te^{-2\int_s^t h(r)\, dr} f(s)\, ds\right)^{1/2}}\, dt =+\infty\, ,\quad \hbox{ for all }\ a\in \R\, ,
\end{equation}
with $h$ defined as above. 
 Secondly, we  remark that     condition \refe{ns-} on the growth of $f$ for $t\to +\infty$ may be written in the equivalent form
\begin{equation}\label{ns-a}
\int_a^{+\infty} \left[ \frac{1}{\left( \int_a^tf(s)\, ds\right)^{1/2}} +\frac{1}{f(t)^{1/q}}\right]\, dt =+\infty\ \hbox{ for all }\ a\in \R\, .
\end{equation} 
Now, we observe that, for $t\geq a$,  on the one hand one has
$$
 \int_a^t e^{-2\int_s^t h(r)\, dr} f(s)\, ds  \leq \int_a^t f(s)\, ds 
$$
and, on the other hand, by the monotonicity of $f$ and $h$,
$$
 \int_a^t e^{-2\int_s^t h(r)\, dr} f(s)\, ds \leq \frac{f(t)}{h(t)} \int_a^t e^{-2 \int_s^t h(r)\, dr} h(s)\, ds\leq \frac{f(t)}{h(t)}=\left( \frac{f(t)}{g^-(t)}\right)^{2/q}\, .
$$
Therefore, for $t\geq a$, one has
\begin{equation}\label{eq1}
\int_a^{+\infty}  \frac{dt}{\left(\int_a^te^{-2\int_s^th(r)\, dr} f(s)\, ds\right)^{1/2}}  \geq \frac{1}{2} \int_a^{+\infty} \left[  \frac{1}{\left( \int_a^tf(s)\, ds\right)^{1/2}}+
 \left( \frac{g^-(t)}{f(t)}\right)^{1/q}\right]\, dt\, .
\end{equation}
 Since $g^-(t)\geq g^-_\infty >0$, it is then clear that
  if \refe{ns-a} is satisfied then \refe{nsh}  holds true, that is to say \refe{ns-} implies \refe{ns0}.

Conversely,  by integration by parts  we also have
$$
\int_a^t f(s)\, ds - \int_a^t e^{-2\int_s^th(r)\, dr} f(s)\, ds = 2\int_a^t   h(s)\left( \int_a^s e^{-2\int_\sigma^sh(r)\, dr} f(\sigma)\, d\sigma \right)\,  ds\, ,
$$
and, by the non increasing monotonicity of the function $h(s)/g^-(s)= \left( \frac{g^-(s)}{f(s)}\right)^{2/q -1}$, it follows that
$$
\int_a^t f(s)\, ds - \int_a^t e^{-2\int_s^th(r)\, dr} f(s)\, ds\leq 2\int_a^t   g^-(s)\left( \int_a^s e^{-2\int_\sigma^sh(r)\, dr} h(\sigma)\, \frac{f(\sigma)}{g^-(\sigma)} \, d\sigma\right)\,  ds\, .
$$
 H\"older inequality with exponent $2/q$ then yields
 $$
\begin{array}{l}
\dis \int_a^t f(s)\, ds -  \int_a^t e^{-2\int_s^th(r)\, dr} f(s)\, ds \\[2ex]
 \dis \leq   \int_a^t g^-(s) \left( \int_a^s 2\, e^{ -2\int_\sigma^s h(r)\, dr}h(\sigma)\left( \frac{f(\sigma)}{g^-(\sigma)}\right)^{2/q} d\sigma \right)^{q/2} \left( \int_a^s 2\, e^{ -2\int_\sigma^s h(r)\, dr}h(\sigma)\, d\sigma \right)^{1-q/2} ds\\[2ex]
 \dis \leq 2^{q/2} \int_a^t g^-(s) \left( \int_a^s e^{ -2\int_\sigma^s h(r)\, dr}f(\sigma)\, d\sigma \right)^{q/2}ds
\end{array}
$$
and then, by monotonicity,
$$
\int_a^t f(s)\, ds -  \int_a^t e^{-2\int_s^th(r)\, dr} f(s)\, ds  \leq 2^{q/2} \left( \int_a^t g^-(s)\, ds\right)\,  \left( \int_a^t e^{ -2\int_\sigma^t h(r)\, dr}f(\sigma)\, d\sigma \right)^{q/2}\, .
$$
By applying further Young inequality with exponent $2/q$ we then obtain
$$
\begin{array}{l}
\dis \int_a^t f(s)\, ds - \int_a^t e^{-2\int_s^th(r)\, dr} f(s)\, ds\\[2ex]
\dis 
   \leq q \,  \int_a^t e^{ -2\int_\sigma^t h(r)\, dr}f(\sigma)\, d\sigma  \frac{\left(\int_a^t g^-(s)\, ds\right)^{2/q}}{\left(\int_a^t f(s)\, ds\right)^{2/q-1}} 
+ \left( 1-\frac{q}{2}\right)  \int_a^t f(s)\, ds\, .
\end{array}
$$
Hence, 
$$
\frac{q}{2}\,   \int_a^t f(s)\, ds \leq \int_a^t e^{-2\int_s^th(r)\, dr} f(s)\, ds \left[ 1+q \frac{\left(\int_a^t g^-(s)\, ds\right)^{2/q}}{\left(\int_a^t f(s)\, ds\right)^{2/q-1}} \right]\, ,
$$
which immediately implies
\begin{equation}\label{eq2}
\int_a^{+\infty} \frac{dt}{\left( \int_a^t e^{-2\int_s^t h(r)\, dr}f(s)\, ds\right)^{1/2}}dt  \leq \frac{2}{\sqrt{q}} \int_a^{+\infty}\left[ \frac{1}{\left( \int_a^t f(s)\, ds\right)^{1/2}} +
\left( \frac{\int_a^t g^-(s)\, ds}{\int_a^t f(s)\, ds}\right)^{1/q}\right]\, dt\, .
\end{equation}
Since for $t$ sufficiently large we have $\int_a^t f(s)\, ds\geq t/2\, f(t/2)$ and, moreover, $\int_a^t g^-(s)\, ds\leq g^-(a) \, t$,  it then easily follows from inequality \refe{eq2} that if \refe{nsh} is satisfied, then \refe{ns-a} must be true.

\hfill$\Box$

\noindent {\bf Proof of Lemma 2.4 } For any fixed $a\in \R$ let us define the real function
$$
\psi_a (t)= \int_a^t e^{-2\int_s^t h(r)\, dr}f(s)\, ds\, ,
$$
which is  the solution of the linear first order initial value problem
$$
\left\{
\begin{array}{c}
\psi_a'(t)+2\, h(t)\, \psi_a(t) =f(t)\\
\psi_a(a)=0
\end{array}\, .
\right.
$$
$\psi_a$ is a $C^1(\R)$ increasing function, since, by the monotonicity of the functions $f$ and $h$, for $t\geq a$, one has
$$
2\, h(t)\, \psi_a(t) \leq 2\, f(t) \int_a^t e^{-2\int_s^t h(r)\, dr}h(s)\, ds < f(t)\, .
$$
Therefore, the condition
$$
{\mathcal R}(a) =\int_a^{+\infty} \frac{dt}{\sqrt{\psi_a(t)}}<+\infty
$$
 only depends on the growth of $\psi_a(t)$ for $t\to +\infty$.   
 We observe that, for any $a\, ,\ a'\in \R$, one has
 \begin{equation}\label{aap}
\psi_{a'}(t) = \psi_a(t) -e^{-2 \int_{a'}^t h(r)\, dr} \psi_a(a')\, ,
\end{equation}
Now, assume there exists $a\in \R$ such that ${\mathcal R}(a)<+\infty$. Then,  it follows that $\lim_{t\to +\infty} \psi_a(t)=+\infty$, and, by the above identity,
$$
\lim_{t\to +\infty} \frac{\psi_{a'}(t)}{\psi_a(t)}=1\, .
$$
Hence, ${\mathcal R} (a')<+\infty$ for all $a'\in \R$, and the first assertion in the statement follows.

Next, let us assume ${\mathcal R}(a) <+\infty$ for all $a\in \R$. For $a_1<a_2$, we easily obtain
$$
{\mathcal R}(a_1) =\int_{a_2}^{+\infty} \frac{d\tau}{\left( \int_{a_2}^\tau e^{-2 \int_\sigma^\tau h(\rho +a_1-a_2)\, d\rho} f(\sigma +a_1-a_2)\, d\sigma \right)^{1/2}} \geq {\mathcal R}(a_2)\, ,
$$
being, by monotonicity, $f(\sigma+a_1-a_2)\leq f(\sigma)$ and $h(\rho +a_1-a_2)\geq h(\rho)\, .$ (We also notice that the above inequality is strict if either $f$ or $h$ is strictly monotone)

Moreover, for any natural number $n$, we have
$$
{\mathcal R}(n)=\int_0^\infty \xi_n (\tau)\, d\tau\, ,
$$
with
$$
\xi_n(\tau) = \frac{1}{\left(\int_0^\tau e^{-2\int_\sigma^\tau h(\rho+n)\,d\rho}\,f(\sigma+n)\,d\sigma\right)^{1/2}}\, .
$$
 $\{ \xi_n\}$ is a non increasing  sequence of integrable functions and, by the monotone convergence theorem, it follows that
$$
\lim_{a \to \infty}{\mathcal R}(a)= \int_0^\infty \frac{d\tau}{\left(f_\infty\int_0^\tau e^{-2 h_\infty (\tau-\sigma)}\,d\sigma\right)^{1/2}}\, ,
$$
with $f_\infty =\lim_{t\to +\infty} f(t)$ and $h_\infty =\lim_{t\to +\infty} h(t)$. If $f_\infty<+ \infty$, then 
$$
\lim_{a \to \infty}{\mathcal R}(a)= +\infty,
$$
which is a contradiction to the non increasing monotonicity of $\mathcal{R}$. Therefore, if ${\mathcal R}(a) < \infty$, then we necessarily have that  $f_\infty=+ \infty$ and \refe{infy} is satisfied.

 \hfill$\Box$

\begin{oss}\label{tutte} {\rm As a byproduct of Theorem \ref{R} we in particular obtain that either all maximal solutions of Cauchy problem \refe{Cpode} blow up for finite $R$, or all maximal solutions are globally defined on $[0,+\infty)$, independently on the initial datum $a$, but only according to the growth for $t\to +\infty$ of $f$ and $g$ measured through conditions \refe{ns+}, \refe{ns0} and \refe{ns-} respectively.}
\end{oss}

\begin{oss}\label{commentns+} {\rm 
As a consequence of the nondecreasing monotonicity of the functions $t\, f(t)$ and $t\, g^+(t)$, condition \refe{ns+} can be easily proved to be equivalent  to the two conditions
$$
 \int_0^{+\infty} \frac{dt}{(tf(t))^{1/2}} =+\infty  \quad \hbox{and }\  \int_{t_0}^{+\infty} \frac{dt}{(t\, g^+(t))^{1/(2-q)}} =+\infty\, ,
$$
for every $t_0>0$ such that $g^+(t_0)>0$. 
If $0< \lim_{t\to +\infty} g(t)<+\infty$, the second integral is infinite because $q\leq1$ and therefore in this case \refe{ns+} becomes the usual Keller--Osserman condition
\begin{equation}\label{KO}
\int_0^{+\infty} \frac{dt}{(t\, f(t))^{1/2}} =+\infty\, .
\end{equation}
If $\lim_{t\to +\infty} g(t)=+\infty$, then \refe{ns+} restricts the growth  at infinity both for  $f$ and $g$; for instance, for power like blowing up functions $g(t)\simeq t^\alpha$ for $t\to+\infty$, then \refe{ns+} requires $\alpha \leq 1-q$, and for $q=1$ at most logarithmic growth $g(t)\simeq (\ln t)^\alpha$ with $\alpha \leq 1$ is allowed.
}
\end{oss}

\begin{oss}\label{commentns-} {\rm On the other hand, 
 condition \refe{ns-} amounts to requiring that either
$$
 \int_0^{+\infty} \frac{dt}{(tf(t))^{1/2}} =+\infty  \quad \hbox{or }\  \int_0^{+\infty} \frac{dt}{(f(t))^{1/q}} =+\infty\, ,
$$
and, for $q=2$,  it is equivalent to
$$
\int_0^{+\infty} \frac{dt}{(f(t))^{1/2}} =+\infty\, .
$$
For functions $f$ having a power growth at infinity, say $f(t)\simeq t^\alpha$ for $t\to +\infty$, \refe{ns-} means that $0\leq \alpha \leq \max 
\{ 1,q \}$, but it includes also functions of the form $f(t)\simeq t\, (\ln t)^\alpha$ with arbitrary $\alpha \geq 0$ if $q>1$, and $0\leq \alpha \leq 2$ if $q\leq 1$.}
\end{oss}

\begin{oss}\label{commentns0} {\rm Let us consider the case $\lim_{t\to +\infty} g(t)=0$. By inequality \refe{eq1},  it follows that 
\begin{equation}\label{s0}
\int_0^{+\infty} \left[ \frac{1}{\left( t\, f(t)\right)^{1/2}} +\left( \frac{g^-(t)}{f(t)}\right)^{1/q}\right]\, dt =+\infty
\end{equation}
is a sufficient condition in order to have that all maximal solutions are globally defined in $[0,+\infty)$. On the other hand, inequality \refe{eq2} shows that if there exists a global maximal solution $\varphi \in C^2([0,+\infty))$, then
\begin{equation}\label{n0}
\int_0^{+\infty} \left[ \frac{1}{\left( t\, f(t)\right)^{1/2}} +\left( \frac{\int_0^tg^-(s)\, ds}{\int_0^t f(s)\, ds}\right)^{1/q}\right]\, dt =+\infty\, .
\end{equation}
Conditions \refe{s0} and \refe{n0} in general are not equivalent  (also for  $q<2$). An easy example occurs for $q=2$, $f(t)\simeq t\, \left( \ln t\right)^\alpha$  and $g(t)\simeq - 1/t$ for $t\to +\infty$. In this case, \refe{s0} is satisfied if and only if $\alpha\leq 2$ whereas \refe{n0} holds true up to $\alpha\leq 3$. We observe that, in this case, \refe{ns0} actually requires $\alpha \leq 2$. On the other hand, if $q<2$ and $g^-(t)$ decays to 0 as a power function of $1/t$, that is $c_1/t^\beta \leq g^-(t)\leq c_2/t^\beta$ for positive constants $c_1,\ c_2$ and $\beta$ and for $t$ sufficiently large, then \refe{s0} and \refe{n0} can be easily proved to be equivalent and, in such a case, they both are a more easy to read equivalent formulation of \refe{ns0}.
}
\end{oss}

\begin{oss}\label{consiKO} {\rm Let us explicitly remark that for $g\equiv 0$, condition \refe{ns0} reduces to the classical Keller--Osserman condition \refe{KO}. An analogous condition is also recovered when $g >0$ in the limit $q\to 0$. Indeed, for $q\to 0$ condition \refe{ns+} becomes \refe{KO} applied to the positive non decreasing nonlinearity $f(t)+g(t)$.
}
\end{oss}

\section{ Viscosity subsolutions of fully nonlinear degenerate elliptic equations}

In this section we apply the previous results on the Cauchy problem \refe{Cpode} to derive a priori  upper estimates and necessary and sufficient conditions for the existence of entire viscosity solutions of inequalities of the form
\begin{equation}\label{eqF}
F(x,D^2u)\geq f(u) +g(u)\, |Du|^q \, ,
\end{equation}
where $F: \R^n\times \mathcal{S}^n\to \R$ is a continuous functions which we assume  always to satisfy  $F(x,O)=0$ and the (normalized) degenerate ellipticity condition
$$
0\leq F(x,X+Y)-F(x,X)\leq \hbox{tr}(Y)\, ,\quad \forall\, x\in \R^n\, ,\ X\, ,\ Y\in \mathcal{S}^n\,,\ Y\geq O\, ,
$$
$\mathcal{S}^n$ being the space of symmetryic $n\times n$ real matrices equipped with the usual ordering.

Special attention will be devoted to study in particular subsolutions of the equation
\begin{equation}\label{eqM+}
\Mpiu (D^2u) = f(u) +g(u)\, |Du|^q
\end{equation}
and  of the equation
\begin{equation}\label{eqPk}
\Pk (D^2u)=  f(u) +g(u)\, |Du|^q\, .
\end{equation}
Let us immediately observe that, by the maximality of operator $\Mpiu$ in the class of second order degenerate elliptic operators,  if $u$ is a viscosity solution of \refe{eqF} then $u$ is a subsolution of \refe{eqM+}.

As a first result, we show that radially symmetric solutions of the above equations can be obtained from  solutions $\varphi \in C^2\left([0,R)\right)$ of problem \refe{Cpode}.

\begin{lemma}\label{L1} (i) Let $f$ and $g$ be continuous non decreasing functions, with $f$ positive.  For any
$q>0$, if $\varphi \in C^2([0,R))$ is a solution of the Cauchy problem \refe{Cpode} with $c=n$, then $\Phi (x)=\varphi (|x|)\in C^2(B_R)$ is a classical solution of equation \refe{eqM+} in the ball $B_R$.

(ii) Let $f$ and $g$ be continuous non decreasing functions, with $f$ positive and  $g$ nonnegative.  For any
$q>0$, if $\varphi \in C^2([0,R))$ is a solution of the Cauchy problem \refe{Cpode} with $c=k$, then $\Phi (x)=\varphi (|x|)\in C^2(B_R)$ is a classical solution of equation \refe{eqPk} in the ball $B_R$.
\end{lemma}
\noindent {\bf Proof}.  A direct computation shows that, if $\Phi (x)=\varphi(|x|)$,  then
$$
D^2\Phi (x)=\left\{
\begin{array}{ll}
\varphi''(0)\, I_n & \hbox{ if } x=0\\[2ex]
\frac{\varphi'(|x|)}{|x|}\, I_n + \left( \varphi''(|x|)-\frac{\varphi'(|x|)}{|x|}\right) \frac{x}{|x|}\otimes \frac{x}{|x|} & \hbox{ if } x\neq 0
\end{array} \right. 
$$
Since $\varphi\in C^2([0,R))$, $\varphi '(0)=0$ and $\varphi''(0)=\lim_{r\to 0} \varphi'(r)/r$, then $\Phi$ belongs to 
 $C^2(B_R)$. Moreover, $\Phi$ is convex since $\varphi$ is convex and increasing, and the very definition of operator $\Mpiu$ yields
 $$
 \Mpiu (D^2\Phi (x)) = \varphi'' (|x|) +(n-1)\, \frac{\varphi'(|x|)}{|x|}\, .
 $$
 Therefore, if $\varphi$ solves \refe{Cpode} with $c=n$, then $\Phi$ solves \refe{eqM+}.
 
 Analogously, we observe that if $\varphi ''(|x|)\geq \frac{\varphi'(|x|)}{|x|}$, then
 $$
 \Pk (D^2\Phi(x))= \varphi'' (|x|) +(k-1)\, \frac{\varphi'(|x|)}{|x|}\, .
 $$
 Therefore, by \refe{sc}, it follows that if $g\geq 0$ and $\varphi$ solves \refe{Cpode} with $c=k$, then $\Phi$ is a solution of \refe{eqPk}.
  
\hfill$\Box$

\noindent In the next result we recall  that form of comparison principle that will be needed in the sequel. It is an immediate consequence of the definition of sub/supersolution when one of the functions to be compared is smooth. For the general regularizing argument needed to compare merely viscosity sub and super solutions we refer to \cite{CIL}.

\begin{prop}\label{comp} Let $f, g$ be continuous functions, with $f$ strictly increasing and $g$ non decreasing, and  let further  $u\in USC(B_R)$ and $\Phi\in C^2(B_R)$ satisfy
$$
F(x, D^2u) -f(u)-g(u)\, |Du|^q\geq 0\geq F(x,D^2\Phi) -f(\Phi) -g(\Phi)\, |D\Phi|^q\quad \hbox{ in } B_R
$$
and
$$\limsup_{|x| \to R^-}\,\left(u(x)-\Phi(x)\right)\le 0\, .$$ 
Then $u(x) \le \Phi(x)$ for all $x \in B_R$\,.  
\end{prop}
\noindent {\bf Proof}. By contradiction, suppose $u-\Phi$ has a positive maximum achieved at some interior point $x_0\in B_R$. By using $\Phi(x)+u(x_0)-\Phi(x_0)$ as test function at $x_0$ in the definition of viscosity subsolution for $u$ it follows that
$$
F(x_0, D^2\Phi(x_0)) \geq  f(u(x_0)) +g(u(x_0))\, |D\Phi(x_0)|^q\, ,
$$
which, by the strict monotonicity of $f$  and the monotonicity of $g$, contradicts the fact that $\Phi$ is a supersolution.
\hfill$\Box$

Our first main result provides  necessary and sufficient conditions for the existence of entire subsolutions of equation \refe{eqM+}.

\begin{theorem}\label{entire}
Let $f,\ g$ be continuous nondecreasing functions, with $f$ positive and strictly increasing.
\begin{itemize}

\item[(i)] If $\dis \lim_{t\to +\infty} g(t)>0$, then there exists $u\in USC(\R^n)$ entire viscosity subsolution of \refe{eqM+} if and only if condition \refe{ns+} is satisfied.

\item[(ii)] If $\dis \lim_{t\to +\infty} g(t)=0$, then there exists $u\in USC(\R^n)$ entire viscosity subsolution of \refe{eqM+} if and only if condition \refe{ns0} is satisfied.

\item[(iii)] If $\dis \lim_{t\to +\infty} g(t)<0$, then there exists $u\in USC(\R^n)$ entire viscosity subsolution of \refe{eqM+} if and only if condition \refe{ns-} is satisfied.
\end{itemize}
\end{theorem}
\noindent {\bf Proof}. (i) If condition \refe{ns+} is satisfied, then any maximal solution of the Cauchy problem \refe{Cpode} with $c=n$, which is globally defined in $[0,+\infty)$ by Theorem \ref{R} (i), gives, by Lemma \ref{L1}, a smooth entire (sub)solution of equation \refe{eqM+}.\\
Conversely, assume by contradiction that there exists $u\in USC(\R^n)$ entire subsolution of \refe{eqM+} and \refe{ns+} does not hold true. 
Let us consider a  maximal solution $\varphi (r)$ of the Cauchy problem \refe{Cpode} with $c=n$ and $a< u(0)$. Then, again by Theorem \ref{R} (i), $\varphi(r)$ blows up for $r=R(a)\in (0,+\infty)$. On the other hand, by Lemma \ref{L1} and Proposition \ref{comp}, the functions $u$ and $\Phi (x)=\varphi(|x|)$ can be compared in $B_{R(a)}$ and we obtain the absurdity
$$
u(0)\leq \Phi (0)=a <u(0)\, .
$$
In the same way, statements (ii) and (iii) follow by Theorem \ref{R} (ii) and (iii) respectively.

\hfill$\Box$

\begin{oss}\label{entireF} {\rm By the maximality of operator $\Mpiu$, conditions \refe{ns+}, \refe{ns0} and \refe{ns-} are necessary conditions for the existence of entire viscosity solutions of inequality \refe{eqF} respectively in the cases $\dis \lim_{t\to +\infty} g(t)>0$, 
$\dis \lim_{t\to +\infty} g(t)=0$ and $\dis \lim_{t\to +\infty} g(t)<0$.}
\end{oss}

The same proof of Theorem \ref{entire} can be applied to subsolutions of equation \refe{eqPk}. But, in this case, the extra assumption $g(t)\geq 0$ is needed in order to have correspondence between solutions of \refe{Cpode} and radially symmetric solutions of \refe{eqPk}.

\begin{theorem}\label{entirePk}
Let $f,\ g$ be continuous nonnegative nondecreasing functions, with $f$ positive and strictly increasing.
There exists an  entire viscosity subsolution $u\in USC(\R^n)$ of equation \refe{eqPk} if and only if condition \refe{ns+} is satisfied.
\end{theorem}

The comparison argument used in  the proof of Theorem \ref{entire}  actually can be applied in order to estimate from above viscosity solutions of inequality \refe{eqF} in any open subset $\Omega \subset \R^n$. If $\partial \Omega\neq \emptyset$, for $x\in \R^n$,  we set
$$
d(x)=\hbox{dist}(x, \partial \Omega)
$$
to denote the distance function from the boundary $\partial \Omega$. In order to state our universal upper bounds we need to invert strictly decreasing functions $\mathcal{R}:\R\to (0,+\infty)$ such that $\dis \lim_{a\to +\infty} \mathcal{R}(a)=0$. For $b>0$, we will denote by $\mathcal{R}^{-1}(b)$ the unique real number $a$ such that $\mathcal{R}(a)=b$ if it exists, and $\mathcal{R}^{-1}(b)=-\infty$ otherwise.

\begin{theorem}\label{upperestimate}
Let $f,\ g$ be continuous nondecreasing functions, with $f$ positive,  strictly increasing and such that
\begin{equation}\label{buf}
\lim_{t\to +\infty} f(t)=+\infty  \quad \hbox{ if }  \lim_{t\to +\infty} g(t)<+\infty\, .
\end{equation}
 Let further $\Omega \subset \R^n$ be an open domain with non empty boundary, and, for $q\in (0,2]$,  let $u\in USC(\Omega)$ be a viscosity solution of \refe{eqF} in $\Omega$.
\begin{itemize}

\item[(i)] Assume that $\dis \lim_{t\to +\infty} g(t)>0$ and  condition \refe{ns+} is not satisfied. Then, pointwisely in $\Omega$, one has
$$
u(x)\leq \max \{ t_0\,,\ \mathcal{R}^{-1}(d(x))\} \, ,
$$
where $t_0=\inf \{ t\in \R\, :\ g(t)\geq 0\}$ and $\mathcal{R} :\R \to (0,+\infty)$ is defined as
\begin{equation}\label{R(a)+}
\mathcal{R} (a) =\left\{
\begin{array}{ll}
\dis 2\, \left(\frac{n}{2-q}\right)^{1/(2-q)} \int_a^{+\infty} \frac{dt}{  \left( \int_a^t f(s)\, ds\right)^{1/2}+\left( \int_a^t g^+ (s)\, ds\right)^{1/(2-q)} } & \hbox{ if } q<2\\[6ex]
\dis  \sqrt{\frac{n}{2}} \int_a^{+\infty} \frac{dt}{ \left(  \int_a^t e^{\frac{2}{n}\int_s^t g^+ (\tau)\, d\tau} f(s)\, ds\right)^{1/2}} & \hbox{ if } q=2
\end{array} \right.\, .
\end{equation}

\item[(ii)] Assume $\dis \lim_{t\to +\infty} g(t)=0$ and  condition \refe{ns0} is not satisfied. Then, pointwisely in $\Omega$, one has
$$
u(x)\leq \mathcal{R}^{-1}(d(x))\, ,
$$
where $\mathcal{R} :\R \to (0,+\infty)$ is defined as
\begin{equation}\label{R(a)0}
\mathcal{R} (a) = \frac{n^{1/q}}{\sqrt{q}}\int_a^{+\infty} \frac{dt}{\left( \int_a^t e^{-2\int_s^t \left(\frac{g^-(r)}{f(r)}\right)^{2/q}f(r)\, dr} f(s)\, ds\right)^{1/2}}\, .
\end{equation}

\item[(iii)] Assume $\dis \lim_{t\to +\infty} g(t)<0$ and  condition \refe{ns-} is not satisfied. Then, pointwisely in $\Omega$, one has
$$
u(x)\leq \mathcal{R}^{-1}(d(x))\, ,
$$
where $\mathcal{R} :\R \to (0,+\infty)$ is defined as
\begin{equation}\label{R(a)-}
\mathcal{R} (a) = 2 \frac{n^{1/q}}{q} \int_a^{+\infty}  \left[ \frac{1}{\left( \int_a^t f(s)\, ds\right)^{1/2}} + \left( \frac{\int_a^tg^-(s)\, ds}{\int_a^t f(s)\, ds}\right)^{1/q}\right]\, dt \, .
\end{equation}

\end{itemize}
\end{theorem}
\noindent {\bf Proof}. (i) As already observed,  $u$ is also a subsolution in $\Omega$ of equation \refe{eqM+}. Let $x_0\in \Omega$ be fixed, and let us set $d_0=d(x_0)$. In particular, $u\in USC\left(\overline{B}_R(x_0)\right)$ is a subsolution of \refe{eqM+} in $B_R(x_0)$ for every $0<R<d_0$. On the other hand, we know from Theorem \ref{R} (i) that, for any $a\in \R$, any maximal solution $\varphi_a$ of \refe{Cpode} is defined on a maximal bounded interval $[0, R(a))$, with $0<R(a)<+\infty$. Moreover, by \refe{blowupfi} and by estimate \refe{ubrg+} of Lemma \ref{stimer}, it follows that, if $a\geq t_0$, then $R(a)\leq \mathcal{R}(a)$, with $\mathcal{R}(a)$ defined by \refe{R(a)+}. Also, Lemma \ref{L1} implies  that the function $\Phi_a(x)=\varphi_a(|x-x_0|)$ solves equation \refe{eqM+} in $B_{R(a)}(x_0)$. Therefore, by Proposition \ref{comp}, it follows that
\begin{equation}\label{upper1}
u(x_0)\leq \Phi_a (x_0)=a\, ,
\end{equation}
provided that $a\geq t_0$ is such that $\mathcal{R}(a)< d_0$. Now, by arguing as in the proof of Lemma \ref{R(a)}, the function $\mathcal{R}$, which is well defined by definition if $q>1$ and by assumption if $q\leq 1$, is easily seen to be strictly decreasing. If $\mathcal{R}(t_0) \leq d_0$, then  the conclusion follows by \refe{upper1}. If $\mathcal{R}(t_0)> d_0$, in order to conclude it is enough to show that $\dis \lim_{a\to +\infty} \mathcal{R}(a)=0$.

If $q=2$, we observe that, for $a$ sufficiently large  such that $g^+(a)>0$, one has
$$
\mathcal{R}(a) \leq \sqrt{\frac{n}{2\, f(a)}} \int_a^{+\infty} \frac{dt}{ \left( e^{\frac{2}{n}g^+(a)\, (t-a)}-1\right)^{1/2}}=\frac{\pi\, n}{2\sqrt{f(a)g^+(a)}}\, ,
$$ 
and the conclusion follows by assumption \refe{buf}.

If $1<q<2$, we note that, again for $a$ such that $g^+(a)>0$ and for every $\delta>0$, one has
$$
\begin{array}{ll}
\dis \mathcal{R}(a)  & \dis \leq 2\, \left(\frac{n}{2-q}\right)^{1/(2-q)} \left[ \int_a^{a+\delta} \frac{dt}{\sqrt{f(a)\, (t-a)}} +\int_{a+\delta}^{+\infty}\frac{dt}{\left( g^+(a)\, (t-a)\right)^{1/(2-q)}}\right]\\[3ex]
& \dis = 2\, \left(\frac{n}{2-q}\right)^{1/(2-q)} \left[ 2\, \sqrt{\frac{\delta}{f(a)}} +\frac{2-q}{q-1} \left( \frac{\delta^{1-q}}{g^+(a)}\right)^{1/(2-q)}\right]\, .
\end{array}
$$
By letting $a\to +\infty$ and by using assumption \refe{buf} we obtain that either
$$
\lim_{a\to +\infty} \mathcal{R}(a) \leq \frac{2}{q-1} \left( \frac{n}{\dis \lim_{a\to +\infty}g(a)}\right)^{1/(2-q)} \left( (2-q)\, \delta\right)^{(1-q)/(2-q)}\, ,
$$
or
$$
\lim_{a\to +\infty} \mathcal{R}(a) \leq 4 \left(\frac{n}{2-q}\right)^{1/(2-q)}\sqrt{\frac{\delta}{\dis \lim_{a\to +\infty} f(a)}}\, ,
$$ 
and the conclusion follows by letting either $\delta \to +\infty$ or $ \delta \to 0$ respectively.

Finally, if $0<q\leq 1$, then, since \refe{ns+} is not satisfied, by Remark \ref{commentns+} we have that either
$$
\int_a^{+\infty} \frac{dt}{\left(\int_a^t f(s)\, ds\right)^{1/2}}<+\infty\quad \hbox{ for every } a\in \R\, ,
$$
or
$$
\int_a^{+\infty} \frac{dt}{\left(t\, g^+(t)\right)^{1/(2-q)}}< +\infty \quad \hbox{ for every } a>0\ \hbox{such that } g^+(a)>0\, .
$$
In the former case, we simply observe that
$$
\mathcal{R}(a)   \leq 2\, \left(\frac{n}{2-q}\right)^{1/(2-q)} \int_a^{+\infty} \frac{dt}{\left(\int_a^t f(s)\, ds\right)^{1/2}}
$$
and that the right hand side goes to 0 as $a\to +\infty$ by applying Lemma \ref{R(a)} with $h\equiv 0$.\newline  
In the latter case, we notice that necessarily one has
$$
0=\lim_{t\to +\infty} \left\{
\begin{array}{ll}
\dis  \frac{t^{1-q}}{g(t)} & \hbox{ if }q<1\\[3ex]
\dis  \frac{\ln t}{g(t)} & \hbox{ if }q=1
 \end{array}\right.\, .
$$
Then, by a similar argument as above,  for every $\delta>0$ sufficiently small, we can estimate, up to the constant factor $2\, \left(\frac{n}{2-q}\right)^{1/(2-q)}$,
$$
 \mathcal{R}(a)   \leq  \left[ \int_a^{a+\delta} \frac{dt}{\sqrt{f(a)\, (t-a)}} 
 + \int_{a+\delta}^{2 a} \frac{dt}{\left( g(a)\, (t-a)\right)^{1/(2-q)}}+\int_{2 a}^{+\infty} \frac{dt}{\left(\frac{t}{2}g\left(\frac{t}{2}\right)\right)^{1/(2-q)}} \right]
$$
 By letting first $a\to +\infty$ and then $\delta\to 0$ we get the conclusion also for every $q\leq 1$.
 
 (ii) We repeat the same comparison argument as above, and conclude by using estimate \refe{ubrg-} of Lemma \ref{stimer} and by applying Lemma \ref{R(a)}.
 
 (iii) We argue again as in (i) and (ii), and we further apply inequality \refe{eq2}. The function $\mathcal{R}$ is easily proved also in this case to be strictly decreasing and to tend to 0 as $a\to +\infty$ by analogous arguments  as in (i).
 
 \hfill$\Box$

\begin{oss}\label{counterexample} {\rm 
The hypothesis \refe{buf} is needed in the case $\dis \lim_{t\to +\infty} g(t)>0$ and $q>1$, whereas it is a consequence of the failure of conditions \refe{ns+}, \refe{ns0} and \refe{ns-} in the other cases. It guarantees that the amplitude $R(a)$ of the maximal interval of existence of a maximal solution of the Cauchy problem \refe{Cpode} satisfies
$$
\lim_{a\to +\infty} R(a)=0\, .
$$
To show it necessity, we can use the same  counterexample exhibited in \cite{P}. Indeed, let us assume $q=2$ and $g\equiv 1$, and for any $a\in \R$ let $\varphi \in C^2([0,R(a)))$ be a maximal solution of
$$
\left\{
\begin{array}{c}
\dis \varphi'' +\frac{n-1}{r}\varphi'= f(\varphi) +(\varphi')^2\\[2ex]
\dis \varphi (0)=a\, ,\ \varphi'(0)=0
\end{array}
\right.\, .
$$
Then, the radial function $\Phi (x)=\varphi (|x|)$ satisfies also
$$
\left\{
\begin{array}{c}
\dis -\Delta \Phi +f(\Phi)+|D\Phi|^2=0\quad \hbox{in } B_{R(a)}\\[2ex]
\dis \Phi =+\infty \quad \hbox{on } \partial B_{R(a)}
\end{array}
\right.\, ,
$$
so that, the function $\Psi(x)= e^{-\Phi (x)}$ solves
$$
\left\{
\begin{array}{c}
\dis -\Delta \Psi=f\left( \ln \frac{1}{\Psi}\right)\, \Psi \quad \hbox{in } B_{R(a)}\\[2ex]
\dis \Psi >0 \quad \hbox{in } B_{R(a)}\, ,\ \Psi =0 \quad \hbox{on } \partial B_{R(a)}
\end{array}
\right.\, .
$$
By the properties of the first eigenvalue of the laplacian, it then follows that 
$$
\lim_{t\to +\infty} f(t)\geq \lambda_1\left(B_{R(a)}\right)\, ,
$$
where $\lambda_1\left(B_{R(a)}\right)$ denotes the first eigenvalue of $-\Delta$ in $B_{R(a)}$. Therefore, if we denote by $\lambda_1$ the first eigenvalue of $-\Delta$ in $B_1$ and if $\dis \lim_{t\to +\infty} f(t)<+\infty$, then we deduce
$$
R(a)\geq \sqrt{\frac{\lambda_1}{\lim_{t\to +\infty} f(t)}}\qquad \forall\, a\in \R\, .
$$
}
\end{oss}

\end{document}